\documentclass[paper=letterpaper, onecolumn, fontsize=11pt]{article}

\usepackage[english]{babel} 
\usepackage{amsmath,amsfonts,amssymb} 
\usepackage{latexsym,amscd,verbatim,alltt,array}
\usepackage{euscript}
\usepackage[affil-it]{authblk}
\usepackage{tikz}
\usepackage{amsmath}
\usepackage{amsfonts}
\usepackage{amsmath}
\numberwithin{equation}{section}
\usepackage{enumitem}
\usetikzlibrary{decorations.pathreplacing,decorations.markings}
\usepackage{multicol}

\makeatletter  
\long\def\@makecaption#1#2{%
  \vskip\abovecaptionskip
  \sbox\@tempboxa{#1 #2}%
  \ifdim \wd\@tempboxa >\hsize
    #1 #2\par
  \else
    \global \@minipagefalse
    \hb@xt@\hsize{\hfil\box\@tempboxa\hfil}%
  \fi
  \vskip\belowcaptionskip}
\makeatother 

\tikzset{ 
  on each segment/.style={
    decorate,
    decoration={
      show path construction,
      moveto code={},
      lineto code={
        \path [#1]
        (\tikzinputsegmentfirst) -- (\tikzinputsegmentlast);
      },
      curveto code={
        \path [#1] (\tikzinputsegmentfirst)
        .. controls
        (\tikzinputsegmentsupporta) and (\tikzinputsegmentsupportb)
        ..
        (\tikzinputsegmentlast);
      },
      closepath code={
        \path [#1]
        (\tikzinputsegmentfirst) -- (\tikzinputsegmentlast);
      },
    },
  },
  mid arrow/.style={postaction={decorate,decoration={
        markings,
        mark=at position .6 with {\arrow[#1]{stealth}}
      }}},
}

\newtheorem{theorem}{Theorem}[section]
\newtheorem{lemma}[theorem]{Lemma}
\newtheorem{proposition}[theorem]{Proposition}

\newtheorem{conjecture}[theorem]{Conjecture}

\newtheorem{definition}[theorem]{Definition}
\newtheorem{remark}[theorem]{Remark}
\newtheorem{example}[theorem]{Example}

\newcommand{\qed}{{\unskip\nobreak\hfil\penalty50\hskip1em\hbox{}\nobreak
   \hfil \ensuremath{\Box}\parfillskip=0pt \par}}

\title{\vspace{-2.5cm} Gorenstein homogeneous subrings of graphs}
\author{\small Lourdes Cruz\thanks{\texttt{lcruzg@math.cinvestav.mx}}\ }
\author{\small Enrique Reyes\thanks{\texttt{ereyes@math.cinvestav.mx}}\ }
\author{\small Jonathan Toledo\thanks{\texttt{jtt@math.cinvestav.mx}}\ }
\affil{\vspace{-0.2cm}\scriptsize Departamento de Matem\'{a}ticas\\ Centro de Investigaci\'{o}n y de Estudios Avanzados del Instituto Polit\'ecnico Nacional\\ Apartado Postal 14--740, Ciudad de M\'{e}xico\\ 07000 M\'exico}
\date{\tiny\ }

\begin{document}
\maketitle
\thispagestyle{empty}

\vspace{-10ex}

\parindent=8mm

\begin{abstract}
\noindent
Let $G=(V,E)$ be a connected simple graph, with $n$ vertices such that $S$ is its homogeneous monomial subring. We prove that if $S$ is normal and Gorenstein, then $G$ is unmixed with cover number $\lceil\frac{n}{2}\rceil$ and $G$ has a strong $\lceil\frac{n}{2}\rceil$-$\tau$-reduction. Furthermore, if $n$ is even, then we show that $G$ is bipartite. Finally, if $S$ is normal and $G$ is unmixed whose cover number is $\lceil\frac{n}{2}\rceil$, we give sufficient conditions for $S$ to be Gorenstein.
\end{abstract}


\section{Introduction}
Let $G=(V(G),E(G))$ be a connected simple graph whose vertex set and edge set are $V(G)=\{x_1,\dotsc,x_n\}$ and $E(G)=\{y_1,\dotsc,y_q\}$, respectively. We assume $E(G)\neq \emptyset$. Let $y=\{x_i,x_j\}$ be an edge of $G$, the \textbf{characteristic vector} of $y$ is the vector in $\{0,1\}^n$ such that its $i$-th entry is $1$, its $j$-th entry is $1$, and the remaining entries are zero. We denote by $v_1,\dotsc,v_q$ the characteristic vector of $y_1,\dotsc,y_q$, respectively. Let $R = K[x_1,\dotsc,x_n]$ be a polynomial ring over a field $K$, the \textbf{homogeneous monomial subring} of $G$ is the ring $S= K[x_1t, \dotsc, x_nt, x^{v_1}t, \dotsc, x^{v_q}t,t] \subset R[t]$,
where $t$ is a new variable. We consider $e_1,\ldots ,e_n$ the canonical vectors in $\mathbb{R}^{n}$. Then, $B=\{ (e_1,1), \dotsc, (e_n,1), (v_1,1), \dotsc ,(v_q,1), e_{n+1}\} \subseteq \mathbb{R}^{n}\times \mathbb{R}$ lies in the hyperplane $x_{n+1} = 1$ where $e_{n+1}=(\underbrace{0,\ldots ,0}_n,1)$. Also, 
$$\mathbb{R}_+B=\Big\{ \sum_{i=1}^q \alpha_i(v_i,1)+\sum_{i=1}^n \beta_i(e_i,1)+\lambda e_{n+1}\mid \alpha_i,\beta_i,\lambda \in \mathbb{R}_+\Big\} \subseteq \mathbb{R}^n\times \mathbb{R}.$$
Thus, $S$ is a standard $K$-algebra, where a monomial $x^at^b$ has degree $b$. We assu\-me $S$ has this grading. Since $(e_1,1),\ldots ,(e_n,1),e_{n+1}$ are linearly independent, $\operatorname{aff}(\mathbb{R}_+B)=\mathbb{R}^{n+1}$. A \textbf{vertex cover} is a subset $\cal C$ of $V(G)$ such that $\mathcal{C} \cap e \neq \emptyset$ for each $e \in E(G)$. The \textbf{cover number} of $G$, $\tau (G)$, is the cardinality of a minimum vertex cover. $G$ is called \textbf{unmixed} if every minimal vertex cover has $\tau (G)$ elements. \\
\noindent
A monomial algebra $\cal A$ is \textbf{Gorenstein} if $\cal A$ is Cohen-Macaulay and its canonical module $\omega_{\cal A}$ is a principal ideal. Hochsther (\cite{Hochster})  proved that if $\cal A$ is normal, then $\cal A$ is Cohen-Macaulay. Hence, if $\cal A$ is normal, then $\cal A$ is Gorenstein if and only if $\omega_{\cal A}$ is principal. 
If $G$ is bipartite, then $S$ is Gorenstein if and only if $G$ 
is unmixed (see~\cite{BRV}). In this paper we prove that if $n$ is even, $S$ is normal and Gorenstein, then $G$ is bipartite. Furthermore, we show that if $S$ is normal and Gorenstein, then $G$ is unmixed, $\tau(G)=\lceil\frac{n}{2}\rceil$ and $G$ has a strong $\lceil\frac{n}{2}\rceil$-$\tau$-reduction (see Definition~\ref{Reduc}). Finally, if $S$ is normal and $G$ is unmixed with $\tau(G)=\lceil\frac{n}{2}\rceil$, we give sufficient conditions for $S$ to be Gorenstein.

\section{Preliminaries}
A subset $F$ of $V(G)$ is a \textbf{stable set} if $y\nsubseteq F$ for each $y\in E(G)$. The cardinality of a maximum stable set is denoted by $\alpha(G)$. $G$ is called {\bf well-covered} if every maximal stable set has $\alpha(G)$ elements. Also, $F$ is a (maximal) stable set if and only if $V(G)\setminus F$ is a (minimal) vertex cover. Hence, $\tau(G)+\alpha(G)=|V(G)|$ and $G$ is unmixed if and only if $G$ is well-covered. $G$ is {\bf very well-covered} if $G$ is well-covered and $\alpha(G)=\frac{n}{2}$ (equivalently, $\tau(G)=\frac{n}{2}$). \\
\noindent
A set of induced subgraphs $G_1,\dotsc,G_s$ of $G$ is a $\tau$-\textbf{reduction} of $G$ if $V(G_1),\dotsc,V(G_s)$ is a partition of $V(G)$ and $\tau(G)=\sum_{i=1}^s \tau(G_i)$. In this case, since $\alpha(G)=n-\tau(G)$ and $n=\sum_{i=1}^n |V(G_i)|$, $\alpha(G)=\sum_{i=1}^n \Big( |V(G_i)| - \tau(G_i)\Big) =\sum_{i=1}^n \alpha(G_i)$.

\begin{lemma}\label{lem02}
If $G$ is unmixed with a $\tau$-reduction $G_1,\dotsc,G_s$, then for each $F$ maximal stable $\alpha(G_i)=|F\cap V(G_i)|$ .
\end{lemma}
\noindent
\begin{proof}
Let $F$ be a maximal stable set. Then, $|F\cap V(G_i)|\leq \alpha(G_i)$. Hence, $\alpha(G)=\sum_{i=1}^s \alpha(G_i)\geq \sum_{i=1}^s |F\cap V(G_i)|=|F|$, since $G_1,\dotsc,G_s$ is a $\tau$-reduction of $G$. But $G$ is well-covered, then $|F|=\alpha(G)$. Therefore $\alpha(G_i)=|F\cap V(G_i)|$.   \qed
\end{proof}

\begin{definition}\rm
An edge $y=\{ x,x^\prime\}$ has the {\it property\/} {\bf (P)} if $\{ z,z^\prime\} \in E(G)$ for each pair of edges $\{ x,z\} ,\{ x^\prime,z^\prime \} \in E(G)$.
\end{definition}

\begin{proposition}{\rm \cite[Theorem 1.2]{Favaron}}\label{prop:unmixed_red}
$G$ is very well-covered if and only if there is a $\tau$-reduction $G_1,\ldots ,G_s$ where each $G_i$ is an edge with the property {\bf (P)}.
\end{proposition}

\begin{proposition}{\rm \cite[Theorem 1.1]{Villa}}\label{prop:verywell}
If $G$ is bipartite and unmixed, then $G$ is very well-covered.
\end{proposition}

\begin{proposition}{\rm \cite[Lemma 14]{Ran-Ves}}\label{prop:reduction}
If $G$ is unmixed, with $\tau(G) = \frac{n+1}{2}$, then there exists a $\tau$-reduction $G_1,\dotsc,G_s$ of $G$ such that $G_i\in E(G)$ for $1\le i\le s-1$ and $G_s$ is a $j$-cycle with $j\in\{3,5,7\}$. 
\end{proposition}

\begin{remark}\rm 
In the previous Proposition $G_1,\ldots ,G_{s-1}$ have the property {\bf (P)}.
\end{remark}
\noindent
\begin{proof}
By contradiction, suppose there are $\{x,z\},\{x^\prime,z^\prime\} \in E(G)$ such that $G_i=\{x,x^\prime\}$ and $\{z,z^\prime\} \notin E(G)$. Then, there is a maximal stable set $F$ such that $\{ z,z^\prime\} \subseteq F$. Hence, $|F\cap V(G_i)|=0$, since $F$ is a stable set. A contradiction, by Lemma~\ref{lem02}, since $\alpha(G_i)=1$. Therefore, $G_i$ has the property {\bf (P)}. \qed
\end{proof}

\begin{proposition}{\rm \cite[Corollary 4.3]{BRV}}\label{thm:BRV}
If $G$ is bipartite, then $S$ is Gorenstein if and only if $G$ is unmixed.
\end{proposition}

\begin{proposition}{\rm \cite[Theorem 6.3.5]{BHer}}\label{CanonMod}
If $S$ is normal, then the canonical module of $S$ is given by
\[
\omega_S = \bigl(\{x^at^b \mid (a,b) \in\mathbb{N}B \cap ({\mathbb R}_+B)^\circ\}\bigr),
\]
where $({\mathbb R}_+B)^\circ$ is the interior of ${\mathbb R}_+B$ relative to $\operatorname{aff}({\mathbb R}_+B)$ (the affine hull of ${\mathbb R}_+B$).
\end{proposition}

\section{Gorenstein homogeneous monomial subrings of \linebreak graphs}
By Theorem 1.1.29 and Proposition 1.1.51 in \cite{Mono}, $\mathbb{R}_+B$ has the unique irreducible representation
$\mathbb{R}_+B=H_{\lambda_1}^+\cap \cdots \cap H_{\lambda_{m_1}}^+$
where $H_{\lambda_i}^+=\{ w\in \mathbb{R}^{n+1}\mid w\cdot \lambda_i\geq 0\}$. Also, by Theorem 1.1.44 in \cite{Mono}, if $F_i=H_{\lambda_i}\cap \mathbb{R}_+B$ and $H_{\lambda_i}=\{ w\in \mathbb{R}^{n+1}\mid w\cdot \lambda_i=0\}$ for $1\leq i\leq m_1$, then $F_1,\ldots ,F_{m_1}$ are the facets of $\mathbb{R}_+B$.

\begin{proposition}\label{prop04-Jun}
$\mathbb{R}_+B=H_{(e_1,0)}^+\cap \cdots \cap H_{(e_n,0)}^+\cap H_{(-\ell_1,1)}^+\cap \cdots \cap H_{(-\ell_m,1)}^+$ where $\ell_1,\ldots ,\ell_m \in \mathbb{R}^n$. Also, $w=(\tilde{w},a)\in (\mathbb{R}_+B)^\circ$ if and only if $\tilde{w}\cdot e_i>0$ for $1\leq i\leq n$ and $w\cdot (-\ell_j,1)>0$ for $1\leq j\leq m$.
\end{proposition}
\noindent
\begin{proof}
We have $\mathbb{R}_+B=H_{\lambda_1}^+\cap \cdots \cap H_{\lambda_{m_1}}^+$. We will prove $\lambda_j=(\tilde{\lambda}_j,0)\in \mathbb{R}^n\times \mathbb{R}$ if and only if $H_{\lambda_j}\in \{H_{(e_1,0)},\ldots ,H_{(e_n,0)}\}$. Assume $\lambda_j=(\tilde{\lambda}_j,0)$, then $\tilde{\lambda}_j\cdot e_i\geq 0$, since $(e_i,1)\in \mathbb{R}_+B\subseteq H_{\lambda_j}^+$ for $1\leq i\leq n$. We take $I=\{ i\mid \tilde{\lambda}_j\cdot e_i=0\}$, then $(e_i,1)\in H_{(\tilde{\lambda}_j,0)}$ if and only if $i\in I$. Furthermore, $(v_k,1)\in H_{(\tilde{\lambda}_j,0)}$ if and only if $y_k=\{ x_{i_1},x_{i_2}\}$ with $i_1,i_2\in I$, since $\tilde{\lambda}_j\cdot e_i\geq 0$. Thus, ${\rm dim}\ F_j=|I|+1$ where $F_j=H_{\lambda_j} \cap \mathbb{R}_+B$, since $e_{n+1}\in H_{\lambda_j}$. But $F_j$ is a facet, then $|I|=n-1$. Hence, $H_{(\tilde{\lambda}_j,0)} \in \{ H_{(e_1,0)},\ldots ,H_{(e_n,0)}\}$.  Now, we prove $H_{(e_k,0)}\in \{ H_{\lambda_1},\ldots ,H_{\lambda_{m_1}}\}$. Since $(e_k,0)\cdot e_{n+1}=0$ and $(e_k,0)\cdot (e_i,1)=0$ for $i\neq k$, we have $e_{n+1},(e_i,0)\in A_k:=\mathbb{R}_+B\cap H_{(e_k,0)}$ for $i\neq k$. Also, $(e_k,0)\cdot (e_k,1)=1$ and $(e_k,0)\cdot (v_i,1)\geq 0$, then $\mathbb{R}_+B\subseteq H_{(e_k,0)}^+$. Hence, $A_k$ is a facet of $\mathbb{R}_+B$, so $H_{(e_k,0)}\in \{ H_{\lambda_1},\ldots ,H_{\lambda_{m_1}}\}$. \\
\noindent
Now, we take $\lambda_j=(\tilde{\lambda}_j,a_j)\in \mathbb{R}^n\times \mathbb{R}$ with $a_j\neq 0$. Since $e_{n+1}\in \mathbb{R}_+B\subseteq H_{\lambda_j}^+$, $a_j=e_{n+1}\cdot \lambda_j>0$. Hence, $H_{(\tilde{\lambda}_j,a_j)}=H_{(\delta_j,1)}$ where $\delta_j=\frac{\tilde{\lambda}_j}{a_j}$. \\
\noindent
Therefore, $\mathbb{R}_+B= H_{(e_1,0)}^+\cap \cdots \cap H_{(e_n,0)}^+ \cap H_{(-\ell_1,1)}^+\cap \cdots \cap H_{(-\ell_m,1)}$. \\
\noindent
Now, by Theorem 1.1.44 in \cite{Mono}, $w=(\tilde{w},a)\in (\mathbb{R}_+B)^\circ$ if and only if $w\cdot (-\ell_j,1)>0$ for $1\leq j\leq m$ and $\tilde{w}\cdot e_i=(\tilde{w},a)\cdot (e_i,0)>0$ for $1\leq i\leq n$.
\qed
\end{proof}

\vspace{2ex}

\noindent
\textbf{Notation}. In this section we take $|a|=a\cdot{\mathbf{1}}= \sum_{i=1}^n a_i $, where ${\bf 1}=(1,\ldots ,1)\in \mathbb{R}^{n}$ and $a=(a_1,\dotsc,a_n)\in \mathbb{R}^n$. Furthermore, if $C$ is a cycle, then $\mathbf{1}_C=\sum_{x_i\in V(C)} e_i$.

\begin{lemma}\label{obs1}
Let $w=(\tilde{w},b)$ be a vector in $\mathbb{N}B$ with $\tilde{w}\in\mathbb{N}^n$ and $b\in\mathbb{N}$. Hence, 
\begin{enumerate}
\item[\rm{1)}] $|\tilde{w}|\le 2b$.
\item[\rm{2)}] If $|\tilde{w}|=2b$, then $w\in \mathbb{N}\big( (v_1,1),\ldots ,(v_q,1)\big)$.
\end{enumerate}
\end{lemma}
\noindent
\begin{proof} 
Since $(\tilde{w},b)\in\mathbb{N}B$, $(\tilde{w},b)=\sum_{i=1}^q \alpha_i(v_i,1) +\sum_{i=1}^n \beta_i(e_i,1) +\lambda e_{n+1}$, where $\alpha_i,\beta_i,\lambda \in \mathbb{N}$. Thus, $b=\sum_{i=1}^q \alpha_i +\sum_{i=1}^n \beta_i +\lambda$. Also, $|\tilde{w}| = \tilde{w}\cdot {\bf 1}= \sum_{i=1}^q \alpha_i(v_i \cdot {\bf 1})$ \linebreak $+\sum_{i=1}^n \beta_i(e_i \cdot {\bf 1}) =2\big( \sum_{i=1}^q \alpha_i \big) +\sum_{i=1}^n \beta_i$. Hence, $2b=|\tilde{w}|+\sum_{i=1}^n \beta_i+2\lambda \ge |\tilde{w}|$. Furthermore, if $|\tilde{w}|=2b$, then $\sum_{i=1}^n \beta_i=0$ and $\lambda =0$. Consequently, $\beta_i=0$ for $1\leq i\leq n$. Therefore, $w\in \mathbb{N}\big( (v_1,1),\ldots ,(v_q,1)\big)$.
\qed
\end{proof}

\begin{lemma}\label{lem12May}
If $w=(\tilde{w},b)=\sum_{i=1}^q \alpha_i(v_i,1)+\sum_{i=1}^n \beta_i(e_i,1)+\lambda e_{n+1}$ with $\alpha_i,\beta_i\in \mathbb{R}_+$, $\lambda >0$ and $\tilde{w}\cdot e_j>0$ for each $1\leq j\leq n$, then $w\in (\mathbb{R}_+B)^\circ$.
\end{lemma}
\noindent
\begin{proof}
We have $(v_i,1),(e_j,1)\in \mathbb{R}_+B\subseteq H_{(-\ell_k,1)}^+$ for $1\leq i\leq q$ and $1\leq j\leq n$, then 
$$\Big( \sum_{i=1}^q \alpha_i(v_i,1)+\sum_{i=1}^n \beta_i(e_i,1)\Big) \cdot (-\ell_k,1)\ge 0.$$ Furthermore, $\lambda e_{n+1}\cdot (-\ell_k,1)=\lambda>0$. Hence, $w\cdot (-\ell_k,1)>0$ for $1\leq k\leq m$. Therefore, by Proposition~\ref{prop04-Jun}, $w\in (\mathbb{R}_+B)^\circ$ since $\tilde{w}\cdot e_i>0$ for $1\leq i\leq n$. 
\qed
\end{proof}

\begin{proposition}\label{prop1}
If $\tau$ is a spanning tree of $G$ and $\tilde{e}_\tau:=\sum_{v_i\in E(\tau)} (v_i,1)+e_{n+1}$, then $\tilde{e}_\tau\in\mathbb{N}B\cap(\mathbb{R}_+B)^\circ$.
\end{proposition}
\noindent
\begin{proof} 
We have $\tilde{e}_\tau\in \mathbb{N}B$. Furthermore, $\tilde{e}_\tau\cdot (e_i,0)>0$ for each $1\leq i\leq n$, since $\tau$ is a spanning tree and $E(G)\neq \emptyset$. Hence, by Lemma~\ref{lem12May}, $\tilde{e}_\tau \in (\mathbb{R}_+B)^\circ$.
\qed
\end{proof}

\begin{remark}\rm\label{remark04-Jun}
If $\tau$ is a spanning tree of $G$, then $|E(\tau)|=n-1$.
\end{remark}

\begin{lemma}\label{lem18Jun}
Assume $\omega_S=(x^\alpha t^\beta)$, then $x^{\tilde{w}} t^a\in \omega_S$ if and only if $(\tilde{w},a)-(\alpha,\beta)\in \mathbb{N}B$.
\end{lemma}
\noindent
\begin{proof}
We have $x^{\tilde{w}}t^a\in \omega_S=(x^\alpha t^\beta)$ if and only if $x^{\tilde{w}}t^a=(x^ut^{a^\prime})(x^\alpha t^\beta)$ with $x^ut^{a^\prime}\in S$. Equivalently, $(\tilde{w},a)-(\alpha,\beta)=(u,a^\prime)\in \mathbb{N}B$.
\qed
\end{proof}

\vspace{1ex}

\noindent
In the following results $\ell_1,\ldots ,\ell_m$ are as in Proposition~\ref{prop04-Jun}.

\begin{proposition}\label{principal-ideal}
If $S$ is normal and $\omega_S$ is principal, then $\omega_S=(x^{\mathbf{1}}t^{\beta})$ where $\beta\leq\lfloor\frac{n}{2}\rfloor+1$.
\end{proposition}
\noindent
\begin{proof} 
By Proposition~\ref{CanonMod}, $\omega_S=(x^\alpha t^{\beta})$ with $(\alpha,\beta)\in\mathbb{N}B \cap (\mathbb{R}_+B)^\circ$, since $\omega_S$ is principal. Also, by Proposition~\ref{prop04-Jun}, $\alpha_i := \alpha \cdot e_i > 0$ for $1\leq i\leq n$. Then, $\alpha_i \geq 1$ since $(\alpha,\beta)\in\mathbb{N}B$. We take $b=\max \{|\ell_1|, |\ell_2|,\dotsc,|\ell_m|,n\}$, then $(\mathbf{1},b)=\sum_{i=1}^n (e_i,1)+(b-n)e_{n+1}\in\mathbb{N}B$. By Proposition~\ref{prop04-Jun}, $(\mathbf{1},b)\in (\mathbb{R}_+B)^\circ$, since $\mathbf{1}\cdot e_i=1>0$ for $1\leq i\leq n$ and $(\mathbf{1},b)\cdot (-\ell_j,1)=-\mathbf{1}\cdot \ell_j+b=-|\ell_j|+b\ge 0$ for $1\leq j\leq m$. Thus, by Proposition~\ref{CanonMod}, $x^{\mathbf{1}}t^b \in \omega_S$. So, by Lemma~\ref{lem18Jun}, $(\mathbf{1},b)-(\alpha ,\beta)\in \mathbb{N}B$. But $\alpha_i\ge 1$, then $\alpha_i=1$. Hence, $\alpha=\mathbf{1}$. Now, if $\tau$ is a spanning tree of $G$, then $\tilde{e}_\tau=(v,n)$ where $v=\sum_{v_i\in E(\tau)} v_i$, since $|E(\tau)|=n-1$ (Remark~\ref{remark04-Jun}). Also, by Propositions~\ref{prop1} and~\ref{CanonMod}, $x^vt^n\in \omega_S=(x^{\mathbf{1}}t^\beta)$. Then, by Lemma~\ref{lem18Jun}, $(v-\mathbf{1},n-\beta)\in {\mathbb N}B$. So, by 1) in Lemma~\ref{obs1}
\[
	2(n-\beta)\ge|v - \mathbf{1}|=|v| - |\mathbf{1}|= 2(n-1) - n = n - 2.
\]
Hence, $\beta\leq\frac{n+2}{2}=\frac{n}{2}+1$. Therefore, $\beta\leq\lfloor\frac{n}{2}\rfloor+1$, since $\beta\in \mathbb{N}$.
\qed
\end{proof}

\begin{lemma}\label{lem21-May_2}
If $\tau$ is a spanning tree of $G$, $e\in E(G)$ and $\tau \cup \{e\}$ has an odd cycle $C$, then the characteristic vectors of the edges of $E(\tau)\cup \{e\}$ are linearly independent.
\end{lemma}
\noindent
\begin{proof}
We can assume $E(C)=\{ y_1,\ldots ,y_k\}$ with $y_i=\{ x_i,x_{i+1}\}$ for $1\leq i\leq k-1$ and $y_k=\{ x_k,x_1\}$. Also, we can suppose $E(\tau)\cup \{ e\}=\{ y_1,\ldots ,y_n\}$, since $|E(\tau)|=n-1$. We will do the proof by induction on $n-k$. If $n-k=0$, then $E(\tau)\cup \{e\}=E(C)$. Thus, $\sum_{i=1}^k (-1)^{i+1}v_i=(e_1+e_2)-(e_2+e_3)+\cdots+(e_k+e_1)=2e_1$. So, $e_1\in \mathbb{R}(v_1,\ldots ,v_k)$. Similarly, $e_i\in \mathbb{R}(v_1,\ldots ,v_k)$ for $1\leq i\leq k$. Hence, $v_1,\ldots ,v_k$ are linearly independent. Now, assume $n-k>0$. Then, there is $x\in V(G)$ such that $deg_{\tau \cup \{e\}}(x)=1$. We can suppose $x=x_n\in y_n$ and $x_n\notin y_j$ for $1\leq j\leq n-1$. Thus, $C\subseteq (\tau \setminus \{ x_n\})\cup \{e\}\subseteq G^\prime:=G\setminus x_n$ and $\tau \setminus \{x_n\}$ is a spanning tree of $G^\prime$. Hence, by induction hypothesis, $v_1,\ldots ,v_{n-1}$ are linearly independent, since $E(\tau \setminus \{x_n\})=E(\tau)\setminus \{ y_n\}$. Therefore, $v_1,\ldots ,v_n$ are linearly independent, since $x_n\in y_n$ and $x_n\notin y_j$ for $1\leq j\leq n-1$.
\qed
\end{proof}

\begin{proposition}\label{prop3}
If $G$ is not bipartite, then $\big( \frac{1}{2},\dotsc,\frac{1}{2}\big) =\frac{1}{2}(\mathbf{1}) \in \{ \ell_1,\ldots ,\ell_m\}$.
\end{proposition}
\noindent
\begin{proof} 
We take $\ell :=\frac{1}{2}(\mathbf{1})\in \mathbb{R}^n$, then $(e_i,1)\cdot (-\ell,1)=\frac{1}{2}$, $(v_j,1)\cdot (-\ell,1)=0$ and $e_{n+1}\cdot (-\ell,1)=1$ for $1\leq i\leq n$ and $1\leq j\leq q$. Hence, $\mathbb{R}_+B\subseteq H_{(-\ell,1)}^+$ and $\{ (v_1,1),\ldots ,(v_q,1)\} \subseteq H_{(-\ell,1)}$. Now, since $G$ is not bipartite, there is an odd cycle $C$ of $G$. We take $e\in E(C)$, then $C-e$ is a path and there is a spanning tree $\tau$ such that $C-e\subseteq \tau$. So, we can assume $E(\tau)\cup \{e\}=\{ y_1,\ldots ,y_n\}$, since $|E(\tau)|=n-1$. Now, if $\mathbf{0}=\sum_{i=1}^n \alpha_i (v_i,1)$, then $\sum_{i=1}^n \alpha_iv_i=0$. But, by Lemma~\ref{lem21-May_2}, $v_1,\ldots ,v_n$ are linearly independent, then $\alpha_1=\alpha_2=\cdots =\alpha_n=0$. Hence, $(v_1,1),\ldots ,(v_n,1)$ are linearly independent in $H_{(-\ell,1)}\cap \mathbb{R}_+B$. Thus, $H_{(-\ell,1)}\cap \mathbb{R}_+B$ is a facet of $\mathbb{R}_+B$. Therefore, by Proposition~\ref{prop04-Jun}, $\ell \in \{ \ell_1,\ldots ,\ell_m\}$.
\qed
\end{proof}

\begin{proposition}\label{IndSetSuppHip}
If $I=\{x_{j_1},\dotsc,x_{j_d}\}$ is a maximal stable set and ${\ell}^\prime = \sum_{i=1}^d e_{j_i}$, then $\ell^\prime \in \{ \ell_1,\ldots ,\ell_m\}$.
\end{proposition}
\noindent
\begin{proof} 
Since $I$ is a stable set, ${\ell}^\prime\cdot v_i \leq 1$ for $1\leq i\leq q$. Then, $(-{\ell}^\prime,1)\cdot (v_i,1)\ge 0$. Also, $(-{\ell}^{\prime},1)\cdot e_{n+1}=1$ and $(-{\ell^\prime},1)\cdot (e_j,1)\ge 0$ for $1\leq j\leq n$. So, $\mathbb{R}_+B\subseteq H_{(-{\ell^\prime},1)}^+$. Now, we can assume $I=\{x_1,\dotsc,x_d\}$, then $\ell^\prime =\sum_{i=1}^d e_i$ and ${\cal C}= V(G) \setminus I = \{x_{d+1}, \dotsc, x_n \}$ is a minimal vertex cover. Thus, there exists $y_{i_1},\dotsc,y_{i_{n-d}} \in E(G)$ such that $y_{i_j} \cap {\cal C} = \{x_{d+j}\}$ for $1\leq j\leq n-d$. Hence, if $w:=(\tilde{w},b)=\sum_{i=1}^d \beta_i(e_i,1)+ \sum_{j=1}^{n-d} \alpha_j(v_{i_j},1)$, then $\tilde{w}\cdot e_{d+j}=\alpha_j$ for $1\leq j\leq n-d$. Assume $w=\mathbf{0}$, then $\alpha_j=\tilde{w}\cdot e_{d+j}=0$. Thus, $\beta_i=0$ for $1\leq i\leq d$. Then, $(e_1,1),\dotsc,(e_d,1),(v_{i_1},1),\dotsc,(v_{i_{n-d}},1)$ are linearly independent. Furthermore, $(e_1,1),\dotsc,(e_d,1),(v_{i_1},1),\dotsc,(v_{i_{n-d}},1)\in H_{(-\ell^\prime,1)}$, since $\ell^\prime =\sum_{i=1}^{d}e_i$ and $|y_{i_j}\cap I|=1$ for $1\leq j\leq n-d$. Hence, $\ell^\prime\in \{ \ell_1,\ldots ,\ell_m\}$. 
\qed
\end{proof}

\vspace{2.5ex}
\noindent
Let $\cal C$ be a minimal vertex cover of $G$. We can suppose ${\cal C}=\{x_1,\dotsc,x_c\}$. Since $\cal C$ is minimal, there exist $y_{r_1},\dotsc,y_{r_c}\in E(G)$ such that $y_{r_i}\cap {\cal C}=\{x_i\}$. We can also suppose $y_{r_i}=\{x_i,x_{j_i}\}$ for $1\leq i \leq c$, where $\{x_{j_1},\dotsc,x_{j_c}\}=\{x_{c+1},\dotsc,x_{c+s}\}$ (some $x_{j_i}$ can be equal to each other). We define
\begin{align}\label{asterisco}
\tilde{e}({\cal C}): &= \sum_{i=1}^c (v_{r_i},1)+\sum_{i=c+s+1}^n (e_i,1)+e_{n+1}
                      = \sum_{i=1}^c (e_i+e_{j_i},1) + \sum_{i=c+s+1}^n (e_i,1)+e_{n+1}\nonumber \\
                     &= \Big( \sum_{i=1}^c e_i + \sum_{i=1}^c e_{j_i} +\hspace{-3mm}\sum_{i=c+s+1}^n\hspace{-3mm}e_i,\hspace{2mm} n-s+1\Big)=\nonumber  \\                   
                     &= (\underbrace{1,\dotsc,1}_c,a_{c+1},\dotsc,a_{c+s},\underbrace{1,\dotsc,1}_{n-(c+s)},n-s+1)
\end{align}
where $a_{c+1},\ldots ,a_{c+s}\ge 1$ and $c=\sum\limits_{i=1}^c \tilde{e} ({\cal C})\cdot (e_i,0)=\sum\limits_{i=1}^c \tilde{e} ({\cal C})\cdot (e_{j_i},0)=\sum\limits_{i=1}^s a_{c+i}$.

\begin{proposition}\label{prop-anterior}
If $\mathcal{C}$ is a minimal vertex cover, then $\tilde{e}({\cal C})\in \mathbb{N}B\cap (\mathbb{R}_+B)^\circ$.
\end{proposition}
\noindent
\begin{proof} 
By~(\ref{asterisco}), $\tilde{e}({\cal C})\cdot (e_i,0)>0$ for $1\leq i\leq n$. Also, $\tilde{e}({\cal C})=\sum_{i=1}^c (v_{r_i},1)+\sum_{j=c+s+1}^n (e_j,1)+e_{n+1}\in\mathbb{N}B$. So, by Lemma~\ref{lem12May}, $\tilde{e}({\cal C})\in(\mathbb{R}_+B)^\circ$.
\qed
\end{proof}

\begin{proposition}\label{thm1}
If $G$ is not bipartite, $S$ is normal and Gorenstein, then $G$ is unmixed, $\tau(G)=\lceil\frac{n}{2}\rceil$ and $\omega_S=(x^\mathbf{1}t^b)$ with $b=\lfloor\frac{n}{2}\rfloor+1$.
\end{proposition}
\noindent
\begin{proof} 
Since $S$ is Gorenstein, $\omega_S$ is principal. So, by Proposition~\ref{principal-ideal}, $\omega_S=(x^\mathbf{1}t^b)$ and $b\le\lfloor\frac{n}{2}\rfloor+1$. Then, by Proposition~\ref{CanonMod}, $(\mathbf{1},b)\in (\mathbb{R}_+B)^\circ$. By Proposition~\ref{prop3}, $\ell =\frac{1}{2}(\mathbf{1})\in \{ \ell_1,\ldots ,\ell_m\}$. Thus, by Proposition~\ref{prop04-Jun}, $(\mathbf{1},b)\cdot(-\ell,1)>0$. So, $b>\ell\cdot\mathbf{1}=\frac{n}{2}$. Hence, $b=\lfloor\frac{n}{2}\rfloor+1$, since $b\in \mathbb{N}$.

\noindent
Now, we prove $G$ is unmixed. Let $\cal C$ a minimal vertex cover. We can assume $\mathcal{C}=\{ x_1,\ldots ,x_c\}$. By Proposition~\ref{prop-anterior}, $\tilde{e}({\cal C})\in \mathbb{N}B\cap (\mathbb{R}_+B)^\circ$. Thus, by Proposition~\ref{CanonMod} and Lemma~\ref{lem18Jun}, $u:=\tilde{e}({\cal C})-(\mathbf{1},b)\in {\mathbb N}B$. Then, by~(\ref{asterisco}), we have 
\[
u = (\underbrace{0,\dotsc,0}_c,a_{c+1}-1,\dotsc,a_{c+s}-1,\underbrace{0,\dotsc,0}_{n-(c+s)},n-s+1-b).
\]   
So, the only possible entries of $u$ different to zero are $c+1,\dotsc,c+s$ and $n+1$. But $\{x_{c+1},\dotsc,x_n\} =V(G)\setminus \mathcal{C}$ is a maximal stable set, then $u=\sum_{i=c+1}^{c+s} \beta_i(e_i,1)+\lambda e_{n+1}$, with $\beta_i,\lambda\in \mathbb{N}$. Then, $u\cdot (\mathbf{1},0)=\sum_{i=c+1}^{c+s} \beta_i$ and $u\cdot e_{n+1}=(\sum_{i=c+1}^{c+s} \beta_i)+\lambda$. Hence, $u\cdot (\mathbf{1},0)\le u\cdot e_{n+1}$ implies $\sum_{i=1}^s (a_{c+i}-1)\le n-s+1-b$. But, by~(\ref{asterisco}), $\sum_{i=1}^s a_{c+i}=c$, then $c-s\le n-s+1-b$. So, $c\leq\lceil\frac{n}{2}\rceil$, since $b=\lfloor\frac{n}{2}\rfloor+1$. Now, by Proposition~\ref{IndSetSuppHip}, $\ell^\prime\in \{ \ell_1,\ldots ,\ell_m\}$ where $\ell^\prime=\sum_{i=c+1}^n e_i$. Thus, by Proposition~\ref{prop04-Jun}, $(-\ell^\prime,1)\cdot(\mathbf{1},b)>0$. Then, $b>\ell^\prime\cdot\mathbf{1}=n-c$. So, $c\ge \lceil\frac{n}{2}\rceil$, since $b=\lfloor\frac{n}{2}\rfloor+1$ and $c\in \mathbb{N}$. Thus, $c=\lceil\frac{n}{2}\rceil$. But $|\mathcal{C}|=c$, then $G$ is unmixed and $\tau(G)=\lceil\frac{n}{2}\rceil$.
\qed
\end{proof}

\begin{proposition}\label{prop04-Jun_cycle}
Assume $S$ is normal and $C$ is an odd cycle with $|V(C)|=k$.
\begin{enumerate}
\item[{\rm 1)}] If $w=(\tilde{w},a)=w^\prime + \big( \mathbf{1}_C,\frac{k+1}{2}\big)$ where $w^\prime \in \mathbb{N}B$, $\tilde{w}\cdot e_i>0$ for each $1\leq i\leq n$, then $x^{\tilde{w}}t^a \in \omega_S$. Furthermore, if $\omega_S=(x^{\mathbf{1}}t^b)$, then $w-(\mathbf{1},b)\in \mathbb{N}B$. 
\item[{\rm 2)}] If $\omega_S=\big( x^{\mathbf{1}}t^\frac{n+1}{2}\big)$ and $\big( \mathbf{1}+e_j,\frac{n+1}{2}\big) \in \mathbb{N}B$, then there is $y\in E(G)$ such that $x_j\in y$ and $y\cap V(C)\neq \emptyset$.
\end{enumerate}
\end{proposition}
\noindent
\begin{proof} 
1) We can assume $C=(y_1,\ldots ,y_k)$ where $\{x_1\}=y_1\cap y_k$. Then, $\big( \mathbf{1}_C,\frac{k+1}{2}\big) =\sum_{\substack{1\le i\le k \\ i \text{ even}}} (v_i,1)+(e_1,1)\in \mathbb{N}B$ (Recall $\mathbf{1}_C=\sum_{x_i\in V(C)} e_i$). Thus, $w\in \mathbb{N}B$ since $w^\prime \in \mathbb{N}B$. Also, $w=w^\prime +\frac{1}{2}\sum_{i=1}^k (v_i,1)+\frac{1}{2}e_{n+1}$, since $\big( \mathbf{1}_C,\frac{k+1}{2}\big) =\frac{1}{2}\sum_{i=1}^k(v_i,1)+\frac{1}{2}e_{n+1}$. Hence, by Lemma~\ref{lem12May}, $w\in (\mathbb{R}_+B)^\circ$ since $\tilde{w}\cdot e_i>0$ for $1\leq i\leq n$. So, by Proposition~\ref{CanonMod}, $x^{\tilde{w}}t^a\in\omega_S$. Now, if $\omega_S=(x^{\mathbf{1}}t^b)$, then by Lemma~\ref{lem18Jun}, $w-(\mathbf{1},b)\in \mathbb{N}B$.  \\
\noindent
2) We take $w_1=(\tilde{w}_1,a_1)=\big( \mathbf{1}+e_j,\frac{n+1}{2}\big)+\big( \mathbf{1}_C,\frac{k+1}{2}\big)$, then $\tilde{w}_1\cdot e_i>0$ for $1\leq i\leq n$. Hence, by 1), $w_1-\big( \mathbf{1},\frac{n+1}{2}\big) \in \mathbb{N}B$, since $\big( \mathbf{1}+e_j,\frac{n+1}{2}\big) \in \mathbb{N}B$. But $w_1-\big( \mathbf{1},\frac{n+1}{2}\big) =\big( \mathbf{1}_C+e_j,\frac{k+1}{2}\big)$. Thus, by 2) in Lemma~\ref{obs1}, $\big( \mathbf{1}_C+e_j,\frac{k+1}{2}\big) =\sum_{i=1}^q \alpha_i(v_i,1)$ with $\alpha_i\in \mathbb{N}$, since $|\mathbf{1}_C+e_j|=2\big( \frac{k+1}{2}\big)$. So, there is $i_1\in \{ 1,\ldots ,q\}$ such that $v_{i_1}=e_j+e_{j^\prime}$ where $x_{j^\prime}\in V(C)$. Therefore, $x_j\in y:=y_{i_1}$ and $y\cap V(C)=\{ x_{j^\prime}\}$. 
\qed
\end{proof}

\begin{theorem} \label{thm2}
If $S$ is normal and $n$ is even, then $S$ is Gorenstein if and only if $G$ is an unmixed bipartite graph.
\end{theorem}
\noindent
\begin{proof} 
$\Rightarrow$) By contradiction suppose $G$ is not bipartite, then $G$ has an odd $k$-cycle $C$. By Proposition~\ref{thm1}, $\omega_S=(x^{\mathbf{1}}t^b)$ where $b=\frac{n}{2}+1$, $\tau (G)=\frac{n}{2}$ and $G$ is unmixed. Then, $G$ is very well-covered and by Proposition~\ref{prop:unmixed_red}, there is a $\tau$-reduction $G_1,\ldots ,G_s$ with $G_i\in E(G)$. We can assume $G_i=y_i$ for $1\leq i\leq s$. We take $w=(\tilde{w},a):=\sum_{i=1}^s (v_i,1)+\big( \mathbf{1}_C,\frac{k+1}{2}\big)$. Since $y_1,\ldots ,y_s$ is a partition of $V(G)$, $\sum_{i=1}^s (v_i,1)=(\mathbf{1},s)$ and $s=\frac{n}{2}$. Thus, $\tilde{w}\cdot e_i>0$ for $1\leq i\leq n$. Hence, by 1) in Proposition~\ref{prop04-Jun_cycle}, $w-(\mathbf{1},b)\in\mathbb{N}B$. But $w-(\mathbf{1},b)=\big( \mathbf{1}_C,\frac{k-1}{2}\big)$, since $b=\frac{n}{2}+1$. Then, by 1) in Lemma~\ref{obs1}, $2\big( \frac{k-1}{2}\big) \geq|\mathbf{1}_C|=|V(C)|=k$. A contradiction, therefore $G$ is bipartite. Also, by Proposition~\ref{thm:BRV}, $G$ is unmixed.
\\
$\Leftarrow$) By Proposition~\ref{thm:BRV}, $S$ is Gorenstein.
\qed
\end{proof}

\begin{definition}\rm
Let $A$ be a subset of $V(G)$. The {\it closed neighbourhood of $A$\/} is $N_G[A]=A\cup \mbox{$\{ x\in V(G)\mid$ there is $y\in E(G)$ such that $y\cap A\neq \emptyset$ and $x\in y\}$}$.
\end{definition}

\begin{definition}\rm\label{Reduc}
A $\tau$-reduction $G_1,\ldots ,G_s$ is a {\it $\lceil\frac{n}{2}\rceil$-$\tau$-reduction\/} if $G_1,\ldots ,G_{s-1}\in E(G)$ and $G_s\in E(G)$ or $G_s\in \{ C_3,C_5,C_7\}$ ($3$-, $5$- or $7$-cycle).

\noindent
A $\lceil\frac{n}{2}\rceil$-$\tau$-reduction $G_1,\ldots ,G_s$ is {\it strong\/} when $G_s\in E(G)$ or if $G_s\in \{ C_3,C_5,C_7\}$, then for each $x\in N_G[G_s]$ and each odd cycle $C$ of $G$, there is an edge $y$ such that $x\in y$ and $y\cap V(C)\neq \emptyset$.
\end{definition}

\begin{proposition}\label{prop2.9_12-May}
If $S$ is normal, $\omega_S=(x^{\mathbf{1}}t^\frac{n+1}{2})$ and $G_1,\ldots ,G_s$ is a $\lceil\frac{n}{2}\rceil$-$\tau$-reduction with $G_s\in \{ C_3,C_5,C_7\}$, then $G_1,\ldots ,G_s$ is strong.
\end{proposition}
\noindent
\begin{proof}
Let $C$ be an odd $k^\prime$-cycle and $x\in N_G[G_s]$. We can suppose $G_1=y_1,\ldots ,G_{s-1}=y_{s-1}$ and $G_s=(y_{j_1},\ldots ,y_{j_k})$. First assume $x\in V(G_s)$, then we can assume $x=x_1\in y_{j_1}\cap y_{j_k}$. Thus, $(\mathbf{1}+e_1,\frac{n+1}{2})=\sum_{i=1}^{s-1} (v_i,1)+  \sum_{\substack{1\le i\le k \\ i \text{ odd}}} (v_{j_i},1)\in \mathbb{N}B$, since $y_1,\ldots ,y_{s-1},V(G_s)$ is a partition of $V(G)$. Now, suppose $x\in N_G[G_s]\setminus V(G_s)$, then there is $y_{j^\prime}=\{x,x^\prime\} \in E(G)$ with $x^\prime \in V(G_s)$. We can suppose $x^\prime \in y_{j_1}\cap y_{j_k}$ and $x=x_1$, then $\big(\mathbf{1}+e_1,\frac{n+1}{2}\big)=\sum_{i=1}^{s-1} (v_i,1)+\sum_{\substack{1\le i\le k \\ i \text{ even}}} (v_{j_i},1)+(v_{j^\prime},1)\in \mathbb{N}B$. \\
\noindent
Hence, in both cases by 2) in Proposition~\ref{prop04-Jun_cycle}, there is $y\in E(G)$ such that $x=x_1\in y$ and $y\cap V(C)\neq \emptyset$. Therefore, $G_1,\ldots ,G_s$ is a strong $\lceil\frac{n}{2}\rceil$-$\tau$-reduction.
\qed
\end{proof}

\begin{theorem}\label{strongred-char}
If $S$ is normal and Gorenstein, then $G$ is unmixed, $\tau(G)=\lceil\frac{n}{2}\rceil$ and $G$ has a strong $\lceil\frac{n}{2}\rceil$-$\tau$-reduction.
\end{theorem}
\noindent
\begin{proof} 
By Propositions~\ref{thm:BRV} and~\ref{thm1}, $G$ is unmixed. Also, by Propositions~\ref{prop:verywell} and~\ref{thm1}, $\tau(G)=\lceil\frac{n}{2}\rceil$. Thus, by Propositions~\ref{prop:unmixed_red} (if $n$ is even) and~\ref{prop:reduction} (if $n$ is odd), there is a $\lceil\frac{n}{2}\rceil$-$\tau$-reduction $G_1,\ldots ,G_s$. If $G_s\in E(G)$, then $G_1,\ldots ,G_s$ is strong. Now, assume $G_s\in \{ C_3,C_5,C_7\}$. So, $G$ is not bipartite and $n$ is odd, since $V(G_1),\ldots ,V(G_s)$ is a partition of $V(G)$. Then, by Proposition~\ref{thm1}, $\omega_S=(x^\mathbf{1}t^\frac{n+1}{2})$. Hence, by Proposition~\ref{prop2.9_12-May}, $G_1,\ldots ,G_s$ is strong. \qed
\end{proof}

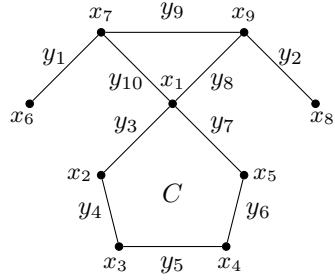
\begin{figure}[h]
\centering
\begin{tikzpicture}[dot/.style={draw,fill,circle,inner sep=1pt},scale=.95]

\node[draw,fill,circle,inner sep=1pt] (1) at (-1,0.5){};
\node (11) at (-1,0.75){{\small $x_{\tiny 7}$}};
\node[draw,fill,circle,inner sep=1pt] (2) at (-2,-.5){};
\node (12) at (-2.1,-.7){{\small $x_{\tiny 6}$}};
\node[draw,fill,circle,inner sep=1pt] (3) at (-1,-1.5){};
\node (13) at (-1.3,-1.5){{\small $x_{\tiny 2}$}};
\node[draw,fill,circle,inner sep=1pt] (4) at (2,-.5){};
\node (14) at (2.1,-.7){{\small $x_{\tiny 8}$}};
\node[draw,fill,circle,inner sep=1pt] (5) at (1,0.5){};
\node (15) at (1,0.75){{\small $x_{\tiny 9}$}};
\node[draw,fill,circle,inner sep=1pt] (6) at (1,-1.5){};
\node (16) at (1.3,-1.5){{\small $x_{\tiny 5}$}};
\node[draw,fill,circle,inner sep=1pt] (7) at (0,-.5){};
\node (17) at (0,-.2){{\small $x_{\tiny 1}$}};

\node[draw,fill,circle,inner sep=1pt] (8) at (-.75,-2.5){};
\node (13) at (-.8,-2.75){{\small $x_{\tiny 3}$}};
\node[draw,fill,circle,inner sep=1pt] (9) at (.75,-2.5){};
\node (16) at (.8,-2.75){{\small $x_{\tiny 4}$}};

\node (y1) at (-1.65,.15){$y_1$};
\node (y2) at (1.65,.15){$y_2$};
\node (y3) at (-0.65,-.8){$y_3$};
\node (y5) at (0,-2.75){$y_5$};
\node (y7) at (0.7,-.8){$y_7$};
\node (y9) at (0,0.8){$y_9$};
\node (y10) at (-0.65,-.2){$y_{10}$};
\node (y8) at (0.7,-.2){$y_8$};
\node (y4) at (-1.15,-2){$y_4$};
\node (y6) at (1.2,-2){$y_6$};
\node (C) at (0,-1.75){$C$};

\path [draw] 
(1) to (2)    (7) to (3)    (8) to (9)
(6) to (7)    (4) to (5)    (8) to (3)
(7) to (5)    (7) to (1)    (5) to (1)
(9) to (6);

\end{tikzpicture}
\caption{$G$ is unmixed and $S$ is not Gorenstein}\label{ExFig}
\end{figure}

\begin{example}\rm\label{Example}
Let $G$ be the graph of the Figure~\ref{ExFig}. Thus, $S$ is normal. Furthermore, $G$ is unmixed with $\tau(G)=\frac{n+1}{2}$ and $y_1,y_2,$ $C=(y_3,y_4,y_5,y_6,y_7)$ is the $\lceil\frac{n}{2}\rceil$-$\tau$-reduction. But it is not strong, since $C^\prime=(y_8,y_9,y_{10})$ is a $3$-cycle and there are not edge between $x_3\in N_G[C]$ and $C^\prime$. Then, by Theorem~\ref{strongred-char}, $S$ is not Gorenstein. Hence, $\omega_S$ is not principal.
\end{example}

\begin{proposition}\label{prop7-May}
Let $y_1,\ldots ,y_u$ be disjoint edges with the property {\bf (P)}. If $w=(\tilde{w},a)=\sum_{i=1}^q \alpha_i(v_i,1)+\sum_{i=1}^n \beta_i(e_i,1)+\lambda e_{n+1}\in \mathbb{N}B$ with $\sum_{i=1}^u \alpha_i$ maximal, then for each $k\in \{1,\ldots ,u\}$, $\alpha_k=\tilde{w}\cdot e_{j_1}$ or $\alpha_k=\tilde{w}\cdot e_{j_2}$, where $y_k=\{ x_{j_1},x_{j_2}\}$. 
\end{proposition}
\noindent
\begin{proof}
Assume $y_k=\{ x_{j_1},x_{j_2}\}$ with $k\in \{ 1,\ldots ,u\}$. First, we prove $\beta_{j_1}=0$ or $\beta_{j_2}=0$. By contradiction, suppose $\beta_{j_1}>0$ and $\beta_{j_2}>0$. Since $(e_{j_1},1)+(e_{j_2},1)=(v_k,1)+e_{n+1}$, we have $w=\sum_{i=1}^q \tilde{\alpha}_i(v_i,1)+\sum_{i=1}^n \tilde{\beta}_i(e_i,1)+\tilde{\lambda} e_{n+1}$ where $\tilde{\alpha}_k=\alpha_k +1$ and $\tilde{\alpha}_i=\alpha_i$ if $i\neq k$; $\tilde{\beta}_{j_1}=\beta_{j_1}-1\ge 0$, $\tilde{\beta}_{j_2}=\beta_{j_2}-1\ge 0$ and $\tilde{\beta}_j=\beta_j$ if $j\notin \{ j_1,j_2\}$; and $\tilde{\lambda}=\lambda +1$. But, $\sum_{i=1}^u \tilde{\alpha}_i =\big( \sum_{i=1}^u \alpha_i\big) +1$. A contradiction, since $\sum_{i=1}^u \alpha_i$ is maximal. Hence, $\beta_{j_1}=0$ or $\beta_{j_2}=0$. We can assume $\beta_{j_2}=0$ and we have two cases: \smallskip

\noindent
\underline{Case $\beta_{j_1}>0$}: We prove $\alpha_l=0$ if $x_{j_2}\in y_l$ and $l\neq k$. By contradiction, suppose there is $y_l\in E(G)$ with $x_{j_2}\in y_l$, $\alpha_l\neq 0$ and $l\neq k$. Then, $l\notin \{1,\ldots ,u\}$, since $y_1,\ldots ,y_u$ are disjoint. We can assume $y_l=\{ x_{j_2},x_{j_3}\}$, then $(v_l,1)+(e_{j_1},1)=(v_k,1)+(e_{j_3},1)$. Thus, $w=\sum_{i=1}^q \alpha^\prime _i(v_i,1)+\sum_{i=1}^n \beta^\prime _i(e_i,1)+\lambda e_{n+1}$ where $\alpha^\prime _k=\alpha_k+1$, $\alpha^\prime _l=\alpha_l-1\ge 0$ and $\alpha^\prime _i=\alpha_i$ if $i\notin \{ k,l\}$; $\beta^\prime _{j_3}=\beta_{j_3}+1$, $\beta^\prime _{j_1}=\beta_{j_1}-1\ge 0$ and $\beta_j^\prime =\beta_j$ if $j\notin \{ j_1,j_3\}$. But, $\sum_{i=1}^u \alpha^\prime _i =\sum_{i=1}^u \alpha_i +1$ since $l\notin \{ 1,\ldots ,u\}$. A contradiction, since $\sum_{i=1}^u \alpha_i$ is maximal. Hence, $\alpha_l=0$ if $x_{j_2}\in y_l$ and $l\neq k$. 

\noindent
\underline{Case $\beta_{j_1}=0$}: We prove $\alpha_l=0$ if $x_{j_1}\in y_l$ with $l\neq k$ or $\alpha_l=0$ if $x_{j_2}\in y_l$ with $l\neq k$. By contradiction, suppose there are $l_1,l_2\in \{ 1,\ldots ,n\} \setminus \{ k\}$ such that $x_{j_1}\in y_{l_1}, x_{j_2}\in y_{l_2},\alpha_{l_1}\neq 0$ and $\alpha_{l_2}\neq 0$. Then, $l_1,l_2\notin \{ 1,\ldots ,u\}$, since $y_1,\ldots ,y_u$ are disjoint. We assume $y_{l_1}=\{ x_{j_1},x_{j^\prime _1}\}$ and $y_{l_2}=\{ x_{j_2},x_{j^\prime _2}\}$. So, $\{ x_{j^\prime _1},x_{j^\prime _2}\}\in E(G)$, since $y_k$ has the property {\bf (P)}. We assume $y_{l^\prime}=\{ x_{j^\prime _1},x_{j^\prime _2}\}$, then $(v_k,1)+(v_{l^\prime},1)=(v_{l_1},1)+(v_{l_2},1)$. Thus, $w=\sum_{i=1}^q \alpha^{\prime \prime}_i(v_i,1)+\sum_{i=1}^n \beta_i(e_i,1)+\lambda e_{n+1}$ where $\alpha^{\prime \prime}_k=\alpha_k+1, \alpha^{\prime \prime}_{l^\prime}=\alpha_{l^\prime}+1, \alpha^{\prime \prime}_{l_1}=\alpha_{l_1}-1\ge 0, \alpha^{\prime \prime}_{l_2}=\alpha_{l_2}-1\ge 0$ and $\alpha^{\prime \prime}_i=\alpha_i$ if $i\notin \{ k,l^\prime,l_1,l_2\}$. But $\sum_{i=1}^u \alpha^{\prime \prime}_i\geq \sum_{i=1}^u \alpha_i+1$, since $l_1,l_2\notin \{ 1,\ldots ,u\}$. A contradiction. Hence, we can suppose $\alpha_l=0$ if $x_{j_2}\in y_l$ and $l\neq k$. \smallskip

\noindent
In both cases, $\tilde{w}\cdot e_{j_2}=\alpha_k$, since $\tilde{w}\cdot e_{j_2}=\sum_{x_{j_2}\in y_l} \alpha_l+\beta_{j_2}$ and $\beta_{j_2}=0$.
\qed
\end{proof}

\begin{definition}\rm\label{Gorens-Reduc}
Let $G_1=y_1,\ldots ,G_{s-1}=y_{s-1},G_s$ be a $\lceil\frac{n}{2}\rceil$-$\tau$-reduction. A representation $w=\sum_{i=1}^q \alpha_i (v_i,1)+\sum_{i=1}^n \beta_i (e_i,1)+\lambda e_{n+1} \in \mathbb{N}B$ (with $\alpha_i,\beta_j,\lambda \in \mathbb{N}$) is {\it principal\/} if it satisfies the following conditions:
\begin{enumerate}
\item[{\rm 1)}] $\sum_{i=1}^u \alpha_i$ is maximal where $u=s$ if $G_s\in E(G)$ or $u=s-1$ if $G_s\notin E(G)$.
\item[{\rm 2)}] If $G_s\in E(G)$, then $\lambda >0$.
\item[{\rm 3)}] If $G_s\in \{ C_3,C_5,C_7\}$, then $G_s=(y_{j_1},\ldots ,y_{j_k})$ such that $\alpha_{j_i}> 0$ for each $i$ even in $\{ 1,\ldots ,k\}$ and $\beta_l>0$ where $x_l\in y_{j_1}\cap y_{j_k}$. 
\end{enumerate}
\end{definition}

\begin{theorem}\label{Gorens}
If $S$ is normal, $G$ is unmixed with a $\lceil\frac{n}{2}\rceil$-$\tau$-reduction and each $w\in (\mathbb{R}_+B)^\circ \cap \mathbb{N}B$ has a principal representation, then $S$ is Gorenstein.  
\end{theorem}
\noindent
\begin{proof}
Let $G_1,\ldots ,G_s$ be a $\lceil\frac{n}{2}\rceil$-$\tau$-reduction. We can assume $G_1=y_1,\ldots ,G_{s-1}=y_{s-1}$; and $G_s=y_s$ if $G_s\in E(G)$. First we prove $x^{\mathbf{1}}t^b\in \omega_S$ with $b=\lfloor \frac{n}{2}\rfloor+1$. If $G_s\in E(G)$, then $n=2s$ and by Lemma~\ref{lem12May}, $(\mathbf{1},b)=\sum_{i=1}^s (v_i,1)+e_{n+1}\in \mathbb{N}B\cap (\mathbb{R}_+B)^\circ$, since $y_1,\ldots ,y_s$ is a partition of $V(G)$. Thus, by Proposition~\ref{CanonMod}, $x^{\mathbf{1}}t^b\in \omega_S$. Now, if $G_s\in \{ C_3,C_5,C_7\}$, then 
\begin{align}\label{asterisco2}
(\mathbf{1},b)=\sum_{i=1}^{s-1} (v_i,1)+\Big( \mathbf{1}_C,\frac{k+1}{2}\Big)
\end{align}
where $C:=G_s$ and $k:=|V(C)|$, since $y_1,\ldots ,y_{s-1},G_s$ is a $\tau$-reduction. Hence, by 1) in Proposition~\ref{prop04-Jun_cycle}, $x^\mathbf{1}t^b\in \omega_S$. \\
\noindent
Now, we take $x^{\tilde{w}}t^a\in \omega_S$ and we prove $x^{\tilde{w}}t^a\in (x^{\mathbf{1}}t^b)$. By Proposition~\ref{CanonMod}, $w:=(\tilde{w},a)\in (\mathbb{R}_+B)^\circ \cap \mathbb{N}B$. Thus, $w$ has a principal representation $w=\sum_{i=1}^q \alpha_i(v_i,1)+\sum_{i=1}^n \beta_i(e_i,1)+\lambda e_{n+1}$. We take $u=s$ if $G_s\in E(G)$ and $u=s-1$ if $G_s\notin E(G)$. So, by Proposition~\ref{prop7-May}, for each $l\in \{ 1,\ldots ,u\}$, $\tilde{w}\cdot e_{i_1}=\alpha_l$ or $\tilde{w}\cdot e_{i_2}=\alpha_l$ where $y_l=\{ x_{i_1},x_{i_2}\}$. Also, by Proposition~\ref{prop04-Jun}, $\tilde{w}\cdot e_{i_1}>0$ and $\tilde{w}\cdot e_{i_2}>0$, since $w\in (\mathbb{R}_+B)^\circ$. Then, $\alpha_l>0$ for $1\leq l\leq u$. If $G_s\in E(G)$, then $u=s$, $\lambda >0$ and $w^\prime:=w-(\mathbf{1},b)=\sum_{i=1}^q \alpha^\prime _i(v_i,1)+\sum_{i=1}^n \beta_i (e_i,1)+\lambda^\prime e_{n+1}$, where $\alpha^\prime _i=\alpha_i-1\geq 0$ if $i\in \{ 1,\ldots ,u\}$, $\alpha^\prime _i=\alpha_i$ in another case; and $\lambda^\prime=\lambda -1\geq 0$. Hence, $w^\prime \in \mathbb{N}B$ implies (by Lemma~\ref{lem18Jun}) $x^{\tilde{w}}t^a\in (x^\mathbf{1}t^b)$. Now, assume $G_s\in \{ C_3,C_5,C_7\}$. Thus, $u=s-1$ and $G_s=C=(y_{j_1},\ldots ,y_{j_k})$ with $\alpha_{j_i}>0$ for each $i$ even in $\{ 1,\ldots ,k\}$ and $\beta_1>0$ where $x_1\in y_{j_1}\cap y_{j_k}$. Since $\big( \mathbf{1}_C,\frac{k+1}{2}\big)=\sum_{\substack{1\le i\le k \\ i \text{ even}}}(v_{j_i},1)+(e_1,1)$, by~(\ref{asterisco2}), $w^{\prime \prime}:=w-(\mathbf{1},b)=\sum_{i=1}^q \alpha^{\prime \prime}_i (v_i,1)+\sum_{i=1}^n \beta^{\prime \prime}_i(e_i,1)+\lambda e_{n+1}$, where $\beta^{\prime \prime}_1=\beta_1-1\geq 0$, $\beta^{\prime \prime}_i=\beta_i$ if $i\neq 1$; $\alpha^{\prime \prime}_i=\alpha_i-1\geq 0$ if $i\in \{ 1,\ldots ,u\}\cup \{ j_2,j_4,\ldots ,j_{k-1}\}$ and $\alpha^{\prime \prime}_i=\alpha_i$ in another case. So, $w^{\prime \prime}\in \mathbb{N}B$. Hence, by Lemma~\ref{lem18Jun}, $x^{\tilde{w}}t^a\in (x^{\mathbf{1}}t^b)$. \\ 
\noindent
Therefore, $\omega_S=(x^\mathbf{1}t^b)$ implies $S$ is Gorenstein, since $S$ is normal.   \qed
\end{proof}

\begin{conjecture}\label{conj_1}
If $S$ is normal and $G$ is unmixed with a strong $\lceil\frac{n}{2}\rceil$-$\tau$-reduction, then each $w\in (\mathbb{R}_+B)^\circ \cap \mathbb{N}B$ has a principal representation.
\end{conjecture}

\begin{conjecture}
Assume $S$ is normal. Hence, $S$ is Gorenstein  if and only if $G$ is unmixed with a strong $\lceil\frac{n}{2}\rceil$-$\tau$-reduction.
\end{conjecture}
\noindent
\begin{proof}
$\Rightarrow$) By Theorem~\ref{strongred-char}. \\
$\Leftarrow$) (Using Conjecture~\ref{conj_1}) By Conjecture~\ref{conj_1}, each $w\in(\mathbb{R}_+B)^\circ \cap \mathbb{N}B$ has a principal representation. Hence, by Theorem~\ref{Gorens}, $S$ is Gorenstein.    \qed
\end{proof}

\bibliographystyle{amsplain}

\end{document}